\documentclass
[11pt,a4paper]
{amsart}
\usepackage[T1]{fontenc}
\usepackage[utf8]{inputenc}
\usepackage{lmodern}
\usepackage{amsmath,amssymb,amsthm,mathtools}
\usepackage{geometry}
\usepackage{microtype}
\usepackage{url}
\numberwithin{equation}{section}
\newtheorem{theorem}{Theorem}[section]
\newtheorem{lemma}[theorem]{Lemma}
\newtheorem{corollary}[theorem]{Corollary}
\theoremstyle{remark}
\newtheorem{remark}[theorem]{Remark}
\theoremstyle{plain}
\newtheorem{proposition}[theorem]{Proposition}
\newtheorem{example}[theorem]{Example}
\newcommand{\BH}{\mathcal{B}(\mathcal{H})}
\newcommand{\HH}{\mathcal{H}}
\newcommand{\KK}{\mathcal{K}}
\newcommand{\EE}{\mathcal{E}}
\newcommand{\ran}{\mathcal{R}}
\newcommand{\ls}{\ell_{+}^{2}}
\begin{document}

\numberwithin{equation}{section}

\title[Higher-Order Isometric Liftings]{Higher-Order Isometric Liftings of \\Convex and Concave Operators}
\author{Laurian Suciu}
\date{}

\begin{abstract}
We study higher-order isometric liftings of bounded operators on Hilbert spaces under convexity and concavity assumptions. First, we prove that every $m$-convex operator $T$ satisfying
$\sup_{n\geq 1}n^{-(m-1)}\|T^n\|^2<\infty$ admits a minimal $m$-isometric lifting whose complementary diagonal entry is an expansive $(m-1)$-isometric weighted shift. Next, we show that every $m$-concave operator admits a minimal $m$-isometric lifting with an expansive $m$-isometric weighted-shift component, without imposing an additional growth assumption. We also establish an alternation between convexity and concavity under polynomial growth conditions and give examples illustrating limitations of converse implications.
\end{abstract}

\keywords{convex and concave operators, m-isometry, lifting}

\subjclass
%%%%%
%%%Mathematics Subject Classification 2020
%%%%%
{47A05, 47A15, 47A20, 47A63}

\dedicatory{Dedicated to Professor L\'aszl\'o K\'erchy
	on the occasion of his 75th birthday,
	with esteem and gratitude}

\maketitle 

\section{Introduction and preliminaries}

\subsection{Introduction}
The study of $m$-isometric operators grew out of Agler's work \cite{Agler} on a disconjugacy theorem for Toeplitz operators and was developed systematically in the three papers of Agler and Stankus \cite{AglerStankusI,AglerStankusII,AglerStankusIII}. An $m$-isometry is characterized by the fact that, for every vector $h$, the sequence $\{\|T^nh\|^2\}_{n\geq0}$ has vanishing $m$th forward difference. This polynomial behavior of the squared norms is one of the basic features of the theory. Weighted shifts provide an important source of examples and structural results; see, for instance, \cite{AbdullahLe,BermudezMartinonNegrin,BuchalaShift}.

A natural question is whether operators satisfying higher-order convexity or concavity conditions can be lifted to $m$-isometries, and whether such liftings can be chosen with useful additional properties. For $m\ge 2$, convex and concave operators have been studied extensively in connection with two-isometric liftings, dilations, and Brownian-type structures; relevant results include \cite{BMS,BS1,BS2,Bu,CS,MajdakBrownianConcave,MajdakBrownianParts,MajdakSuciuConvex,MajdakSuciuTriangulations,SuciuTwoLiftings,SuciuTotoi}. More generally, the existence and structure of expansive $m$-isometric liftings were considered in \cite{SuciuExpansive}. These works provide the background for the present constructions.

Polynomial growth of the powers is also closely connected with higher-order isometric liftings. In \cite{BMS} it was proved that polynomial growth characterizes the existence of an $m$-isometric lifting for some order, and they also studied two-isometric liftings of convex operators. The first result of this paper is a fixed-order statement: if $T$ is $m$-convex and satisfies $\|T^n\|^2=O(n^{m-1})$, then it has a minimal $m$-isometric lifting with a specific triangular representation. The complementary part is an expansive $(m-1)$-isometric weighted shift, and the lifting satisfies the additional invariance property $S^{*n}S^n\mathcal H\subseteq\mathcal H$, for $n\ge 1$. The construction also preserves expansivity from $T$ to $S$. In the case $m=2$, the existence of two-isometric liftings for convex operators is already part of the earlier theory in \cite{BMS}; our construction gives the stated weighted-shift realization.

For concave operators, two-isometric liftings and related Cauchy-dual and Brownian-type questions were studied in \cite{BS1,BS2,MajdakBrownianConcave,MajdakBrownianParts}. For higher orders, Bucha\l{}a recently obtained in \cite{Bu} a minimal $m$-isometric dilations for expansive $m$-concave operators. Our second main result applies to every $m$-concave operator, without assuming expansivity or imposing a separate growth condition. It gives a minimal $m$-isometric lifting whose complementary part is an expansive $m$-isometric weighted shift. The polynomial growth bound for the original operator follows, in turn, from the existence of this lifting. The construction uses positive operator-valued polynomial sequences; the general structure of $m$-isometric shifts with operator weights has also been studied recently in \cite{BuchalaShift}.

Finally, we combine the two lifting theorems with a finite-difference argument to obtain an alternation between convexity and concavity under suitable polynomial growth assumptions. The resulting conditions and lifting conclusions are summarized between convexity and concavity in a corollary. Two weighted-shift examples show that some implications need not hold.

The paper is organized as follows. In Section~2, after a scalar estimate for sequences with non-negative higher differences, we construct minimal $m$-isometric liftings for $m$-convex operators satisfying the natural polynomial growth condition. We also discuss the cases $m=2$ and $m=3$ and give an example in which the complementary part of the lift can be chosen to be an isometry. 

Section~3 is devoted to $m$-concave operators. We construct some minimal $m$-isometric liftings, in the case $m\ge 2$. This result (Theorem \ref{thm:3.2}) is an effective generalization of the classical theorem of Sz.-Nagy-Foias which asserts that every contraction (i.e., operator 1-concave) admits an (1-) isometric lifting; see \cite{SFb}. Also we establish the alternation principle between convexity and concavity, and provide examples, which show that these concepts cannot be inferred one from the other, or transferred from an order to the other, in general. 

\subsection{Preliminaries}
Throughout the paper, all Hilbert spaces are complex. We write $\mathcal B(\mathcal H,\mathcal K)$ for the space of bounded linear operators from a Hilbert space $\mathcal H$ into a Hilbert space $\mathcal K$, and $\mathcal B(\mathcal H)=\mathcal B(\mathcal H,\mathcal H)$. The identity operator on $\mathcal H$ is denoted by $I_{\mathcal H}$, or simply by $I$ when the space is clear. For $A\in\mathcal B(\mathcal H,\mathcal K)$, we denote its range and null space by $\mathcal R(A)$ and $\mathcal N(A)$, respectively. The closure of a set is indicated by an overline. If $A\in\mathcal B(\mathcal H)$ is positive, we write $A\geq0$ and denote its positive square root by $A^{1/2}$. We say that $T\in\mathcal B(\mathcal H)$ is {\it expansive} if $T^*T\geq I_{\mathcal H}$.

For a Hilbert space $\mathcal E$, let
\[
\ell^2_+(\mathcal E)=\left\{(x_0,x_1,\ldots):x_j\in\mathcal E,\ \sum_{j=0}^{\infty}\|x_j\|^2<\infty\right\}.
\]
For each $q\geq0$, $J_q:\mathcal E\to\ell^2_+(\mathcal E)$ denotes the canonical embedding into the $q$th coordinate. A unilateral weighted shift with operator weights $W_j\in\mathcal B(\mathcal E)$ is the operator on $\ell^2_+(\mathcal E)$ given by
\[
W(x_0,x_1,\ldots)=(0,W_0x_0,W_1x_1,\ldots), \quad x_j \in \mathcal{E}, \quad j\ge 0.
\]

For integers $n,p\geq0$, we use the convention $\binom{n}{p}=0$ if $n<p$, and otherwise the usual binomial coefficient. In particular, for fixed $p$, the expression $\binom{n}{p}$ is the restriction to nonnegative integers of a polynomial in $n$ of degree $p$.

For $T\in\mathcal B(\mathcal H)$, we set $\Delta_0(T)=I_{\mathcal H}$ and, for every integer $m\geq1$, define the operator
\[
\Delta_m(T)=\sum_{j=0}^{m}(-1)^{m-j}\binom{m}{j}T^{*j}T^j.
\]
These finite difference operators satisfy
\[
\Delta_{m+1}(T)=T^*\Delta_m(T)T-\Delta_m(T),\qquad m\geq0.
\]
We call $T$ an $m$-{\it isometry} if $\Delta_m(T)=0$, $m$-{\it convex} if $\Delta_m(T)\geq0$, and $m$-{\it concave} if $\Delta_m(T)\leq0$.

For a scalar sequence $a=(a_n)_{n\geq0}$, its forward difference is $\delta a_n=a_{n+1}-a_n$. In particular, if $a_n=\|T^nh\|^2$ for $h\in\mathcal H$, then
\[
(\delta^m a)_n=(T^{*n}\Delta_m(T)T^nh,h),\qquad n,m\geq0.
\]

If $\mathcal H$ is identified with a closed subspace of a Hilbert space $\mathcal K$, an operator $S\in\mathcal B(\mathcal K)$ is called a {\it lifting} of $T\in\mathcal B(\mathcal H)$ if $\mathcal H$ is invariant for $S^*$ and $P_{\mathcal H}S|_{\mathcal H}=T$, where $P_{\mathcal H}$ is the orthogonal projection onto $\mathcal H$. Equivalently, relative to $\mathcal K=\mathcal H\oplus\mathcal L$, such a lifting has a triangular representation
\[
S=\begin{pmatrix}T&0\\E&W\end{pmatrix}
\]
for suitable bounded operators $E:\mathcal H\to\mathcal L$ and $W\in\mathcal B(\mathcal L)$. The lifting is {\it minimal} if
\[
\mathcal K=\bigvee_{n\geq0}S^n\mathcal H,
\]
where the right-hand side denotes the closed linear span. In this paper we also pay particular attention to liftings satisfying $S^{*n}S^n\mathcal H\subseteq\mathcal H$ for every $n\geq1$.

% The remaining text is the transcription of the handwritten manuscript,
% with the corrections explicitly agreed upon in the preceding discussion.

\setcounter{section}{1}
\section{$m$-Isometric Liftings for $m$-Convex Operators}

In this section we prove that an $m$-convex operator for an integer $m\geq2$, which satisfies a necessary boundedness condition, admits a special $m$-isometric lifting. This result is related to the polynomial-growth lifting theorem and the two-isometric lifting results for convex operators in \cite{BMS}, as well as to the study of expansive $m$-isometric liftings in \cite{SuciuExpansive}. 

We begin with a useful result concerning an estimate for the scalar finite differences of higher order. We use the operator $\delta$ defined on sequences $x=\{x_n\}_{n\geq0}$ of elements in a linear space by $(\delta x)_n=x_{n+1}-x_n$.

\begin{lemma}\label{lem:scalar}
Let $a=\{a_n\}_{n\geq0}$ be a sequence of scalars satisfying the conditions
\begin{equation}\label{eq:2.1}
0\leq a_n\leq\rho(n+1)^{m-1},\qquad (\delta^m a)_n\geq0\qquad(n\geq0)
\end{equation}
for a constant $\rho>0$ and an integer $m\geq2$. Then
\begin{equation}\label{eq:2.2}
(\delta^{m-1}a)_n\leq(m-1)!\rho\qquad(n\geq0).
\end{equation}
\end{lemma}
\begin{proof}
From the identities
\[
a_{n+p}-a_n=\sum_{i=0}^{p-1}\delta a_{n+i},\qquad
\sum_{q=j}^{N}\binom{q}{j}=\binom{N+1}{j+1},
\]
where $n,j\geq0$, $p\geq1$ and $N\geq j$, we infer the relation
\[
a_{n+p}=\sum_{j=0}^{m-1}\binom{p}{j}\delta^j a_n
+\sum_{q=0}^{p-m}\binom{p-q-1}{m-1}\delta^m a_{n+q}
\qquad(p\geq m,\ n\geq0).
\]
Since $\delta^m a_{n+j}\geq0$ we obtain for $n\geq0$ and $p\geq1$,
\begin{equation}\label{eq:2.3}
a_{n+p}\geq\sum_{j=0}^{m-2}\binom{p}{j}\delta^j a_n
+\binom{p}{m-1}\delta^{m-1}a_n.
\end{equation}
To estimate $p^{-(m-1)}a_{n+p}$ notice first that
\[
\frac{1}{p^{m-1}}\sum_{j=0}^{m-2}\binom{p}{j}\delta^j a_n
\xrightarrow[p\to\infty]{}0,\qquad
\frac{1}{p^{m-1}}\binom{p}{m-1}
\xrightarrow[p\to\infty]{}\frac{1}{(m-1)!}.
\]
On the other hand, we have from (\ref{eq:2.1}) that
\[
\frac{a_{n+p}}{p^{m-1}}
\leq\rho\left(\frac{n+p+1}{p}\right)^{m-1}
\xrightarrow[p\to\infty]{}\rho.
\]
Finally, by using these limits we obtain by (\ref{eq:2.3}) the inequality (\ref{eq:2.2}).
\end{proof}

Now we prove the main result of this section.
\begin{theorem}\label{thm:main}
Let $T\in\BH$ be $m$-convex for an integer $m\geq2$ and satisfying the condition
\begin{equation}\label{eq:2.4}
\sup_{n\geq1}\frac{1}{n^{m-1}}\|T^n\|^2<+\infty.
\end{equation}
Then $T$ admits a minimal $m$-isometric lifting $S$ on a Hilbert space $\KK=\HH\oplus\mathcal{L}$ of the form
\begin{equation}\label{eq:2.5}
S=\begin{pmatrix}T&0\\E&W\end{pmatrix},
\end{equation}
where $W$ is an expansive $(m-1)$-isometric weighted shift on $\mathcal{L}=\ls(\overline{\ran(E)})$, while $E:\HH\to\mathcal{L}$ is an operator with
\begin{equation}\label{eq:2.6}
E^*Eh=\lim_{n\to\infty}\sum_{j=0}^{n}T^{*j}\Delta_m(T)T^jh\qquad(h\in\HH).
\end{equation}
Furthermore, one has $S^{*n}S^n\HH\subset\HH$ for $n\geq1$, and if $T$ is expansive then $S$ is expansive, too.
\end{theorem}
\begin{proof}
We will do the construction of $S$ in (\ref{eq:2.5}) in three steps.

\emph{Step 1.} We show that the sequence from (\ref{eq:2.6}) strongly converges to a positive operator in $\BH$, which leads to the operator $E$. Indeed, since $\Delta_m(T)=T^*\Delta_{m-1}(T)T-\Delta_{m-1}(T)$ we obtain for integers $n\geq1$ that
\begin{equation}\label{eq:2.7}
T^{*n}\Delta_{m-1}(T)T^n=\Delta_{m-1}(T)
+\sum_{j=0}^{n-1}T^{*j}\Delta_m(T)T^j.
\end{equation}
Put $s_n(T):=\sum_{j=0}^{n-1}T^{*j}\Delta_m(T)T^j\geq0$ because $\Delta_m(T)\geq0$, and so we have $s_n(T)\leq s_{n+1}(T)$. Hence the sequence $\{s_n(T)\}_{n\geq1}$ is increasing.

On the other hand, the condition (\ref{eq:2.4}) shows that there exists a constant $\xi>0$ such that $\|T^n\|^2\leq\xi(n+1)^{m-1}$ for $n\geq0$. Then for each $h\in\HH$, the sequence $a_h=\{\|T^nh\|^2\}_{n\geq0}$ satisfies
\[
(\delta^j a_h)_n=(T^{*n}\Delta_j(T)T^nh,h)
\qquad(1\leq j\leq m,\ n\geq0),
\]
where $(\delta^m a_h)_n\geq0$ for $n\geq0$. Since we also have
\[
a_{h,n}=\|T^nh\|^2\leq\xi\|h\|^2(n+1)^{m-1}
\qquad(n\geq0),
\]
from Lemma \ref{lem:scalar} one obtains
\[
(T^{*n}\Delta_{m-1}(T)T^nh,h)\leq\xi\|h\|^2(m-1)!
\qquad(n\geq0).
\]

This implies by (\ref{eq:2.7}) that the sequence $\{s_n(T)\}_{n\geq0}$ is bounded. Because it is also increasing, it follows that this sequence strongly converges in $\BH$. We denote
\begin{equation}\label{eq:2.8}
D:=\operatorname{s-lim}_{n\to\infty}s_n(T)
=\sum_{j=0}^{\infty}T^{*j}\Delta_m(T)T^j.
\end{equation}
Clearly, we have $0\leq\Delta_m(T)\leq D$ and since $T^*s_n(T)T=s_{n+1}(T)-\Delta_m(T)$ for $n\geq1$, by passing to limit for $n\to\infty$ one obtains
\begin{equation}\label{eq:2.9}
D-T^*DT=\Delta_m(T).
\end{equation}
On the other hand, the above inequality between $\Delta_m(T)$ and $D$ shows us that there exists a contraction $C:\overline{\ran(D)}\to\overline{\ran(\Delta_m(T))}$ such that
\[
CD^{1/2}=[\Delta_m(T)]^{1/2}.
\]

Let $B=C^*C\in\mathcal{B}(\overline{\ran(D)})$ and $E_0:\HH\to\overline{\ran(D)}$ be defined by $E_0h=D^{1/2}h$, for $h\in\HH$. It is clear that
\begin{equation}\label{eq:2.10}
E_0^*E_0=D,\qquad E_0^*BE_0=\Delta_m(T),\qquad 0\leq B\leq I_{\EE},
\end{equation}
where $\EE=\overline{\ran(E_0)}=\overline{\ran(D)}$.

Finally, we define the operator $E:\HH\to\mathcal{L}:=\ls(\EE)$ as $E=JE_0$, where $J:\EE\to\mathcal{L}$ is the embedding mapping. Then $E^*E=D$, so $E$ satisfies the relation (\ref{eq:2.6}).

\emph{Step 2.} At this stage we define the shift $W$, in the case $m\geq3$. We use the operators $B_n\in\mathcal{B}(\EE)$ for $n\geq0$, given by
\begin{equation}\label{eq:2.11}
B_n=(1+s_n)I_{\EE}+nB,\qquad
s_n=\begin{cases}
0,&0\leq n<m-2,\\[2pt]
\displaystyle\binom{n}{m-2},&n\geq m-2.
\end{cases}
\end{equation}
Thus $B_0=I_{\EE}\leq B_n$ for $n\geq1$, which ensures that $B_n$ is invertible in $\mathcal{B}(\EE)$. Now the weighted shift $W$ on $\mathcal{L}=\ls(\EE)$ is defined as
\[
W(e_0,e_1,\ldots)=(0,W_1e_0,W_2e_1,\ldots)
\qquad(e_n\in\EE,\ n\geq0),
\]
where the weights $W_n\in\mathcal{B}(\EE)$ are given by
\[
W_n=B_n^{1/2}B_{n-1}^{-1/2}\qquad(n\geq1).
\]

We want to show that the sequence $\{W_n\}_{n\geq1}$ is bounded. Indeed, we first infer from (\ref{eq:2.11}) that
\begin{equation}\label{eq:2.12}
B_n-B_{n-1}=B+(s_n-s_{n-1})I\qquad(n\geq1).
\end{equation}
But from the estimate of $s_n$ in (\ref{eq:2.11}) we see that $s_n=s_{n-1}=0$ for $n<m-2$, and $s_n-s_{n-1}=1=\binom{n}{m-2}$ when $n=m-2$. In the case $n>m-2$ we have
\[
\begin{aligned}
s_n-s_{n-1}
&=\binom{n}{m-2}-\binom{n-1}{m-2}
=\binom{n-1}{m-3}\\
&=\frac{m-2}{n-m+2}\binom{n-1}{m-2}
\leq(m-2)\binom{n-1}{m-2}=(m-2)s_{n-1}.
\end{aligned}
\]
Thus, for every integer $n\geq1$ one has
\[
s_n-s_{n-1}\leq(m-2)(1+s_{n-1}).
\]
Next, since $B\leq I_{\EE}\le B_{n-1}$ and $B_{n-1}\geq(1+s_{n-1})I_{\EE}$, it follows from (\ref{eq:2.12}) and the above inequality that
\[
B_n\leq B_{n-1}+B +(m-2)(1+s_{n-1})I_{\EE}
\leq mB_{n-1}\qquad(n\geq2).
\]
But $B_1\leq(2+s_1)I_{\EE}\leq mI_{\EE}=mB_0$, since $s_1=1$ if $m=3$ and $s_1=0$ if $m\geq4$. Therefore for every integer $n\geq1$ we have $I_{\EE}\leq B_n\leq mB_{n-1}$,
which gives
\[
I_{\EE}\leq W_n^*W_n
=B_{n-1}^{-1/2}B_nB_{n-1}^{-1/2}\leq mI_{\EE}.
\]
We obtained that the sequence $\{W_n\}_{n\geq1}$ is bounded, hence $W$ is bounded and expansive with $\|W\|^2\leq m$.

At the end of this stage we prove that $\Delta_{m-1}(W)=0$. Indeed, since the weights of $W^n$ are the operators
\[
W_{j+n}\cdots W_{j+1}=B_{j+n}^{1/2}B_j^{-1/2}
\qquad(j,n\geq0),
\]
it follows that the $j$-coordinate of $W^{*n}W^n$ is the operator
\begin{equation}\label{eq:2.13}
(W^{*n}W^n)_j=B_j^{-1/2}B_{j+n}B_j^{-1/2}=:P_j(n).
\end{equation}
But $P_j(n)$ is an operator polynomial in variable $n$ of degree at most $m-2$ (using (\ref{eq:2.10}) for $B_{j+n}$), therefore
\[
\sum_{n=0}^{m-1}(-1)^{m-n-1}\binom{m-1}{n}P_j(n)=0
\qquad(j\geq0).
\]
This then gives $\Delta_{m-1}(W)=0$, that is $W$ is an $(m-1)$-isometry, in the case $m\geq3$.

For $m=2$ the construction of operators $D$ and $E$ is obtained (as before) by means of $\Delta_2(T)\geq0$. In this case $W$ is the forward shift (i.e., $W_n=I_{\EE}$, $n\geq1$) on $\mathcal{L}=\ls(\EE)$, therefore $W$ is a non-unitary isometry with $W^*E=0$ (as $\EE=\overline{\ran(E)}$).

\emph{Step 3.} At this stage we analyze the operator $S$ with the representation (\ref{eq:2.5}) on $\KK=\HH\oplus\ls(\EE)$, where $\EE=\mathcal{N}(W^*)=\overline{\ran(E)}$, so $W^*E=0$.

Assume firstly that $m\geq3$. Since $B_0=I_{\EE}$ we have from (\ref{eq:2.13}) that $J^*W^{*n}W^nJ=B_n$ for $n\geq1$, and this shows that
\[
W^{*n}W^n\mathcal{N}(W^*)\subset\mathcal{N}(W^*)\qquad(n\geq1).
\]
Now by means of (\ref{eq:2.5}) we get for $n\geq1$,
\[
S^n=\begin{pmatrix}T^n&0\\ E_n&W^n\end{pmatrix},
\qquad E_n=\sum_{j=0}^{n-1}W^{n-j-1}ET^j,
\]
such that $W^{*n}E_n=0$ by the above inclusion. We infer that
\begin{equation}\label{eq:2.14}
S^{*n}S^n=\begin{pmatrix}D_n&0\\0&W^{*n}W^n\end{pmatrix},
\qquad D_n=T^{*n}T^n+E_n^*E_n\qquad(n\geq1),
\end{equation}
while $D_0=I_{\HH}$. So, $S^{*n}S^n\HH\subset\HH$ for $n\ge 1$, and $S$ is expansive if $T$ is such.

It can be demonstrated (see lemma below) that $D_n$ has the following expression:
\begin{equation}\label{eq:2.15}
D_n=\sum_{j=0}^{m-1}\binom{n}{j}\Delta_j(T)
+\left[n+\binom{n}{m-1}\right]D\qquad(n\geq1).
\end{equation}
Here $Q_p(n):=\binom{n}{p}$ is a polynomial in $n$ of degree most $p$, where $Q_p(p)=1$ and $Q_p(n)=0$ if $n<p$. Thus, $D_n$ as an operator-valued polynomial in $n$ is of degree at most $m-1$. Hence
\[
\begin{aligned}
\sum_{n=0}^{m}(-1)^{m-n}\binom{m}{n}D_n
&=\sum_{j=0}^{m-1}\left[\sum_{n=0}^{m}(-1)^{m-n}\binom{m}{n}Q_j(n)\right]\Delta_j(T)\\
&\quad+\sum_{n=0}^{m}(-1)^{m-n}\binom{m}{n}[n+Q_{m-1}(n)]D=0,
\end{aligned}
\]
having in view that the difference of order $m$ of each polynomial of degree at most $m-1$ is zero.

Finally, we infer by using (\ref{eq:2.14}) that
\[
\Delta_m(S)=\left[\sum_{n=0}^{m}(-1)^{m-n}\binom{m}{n}D_n\right]
\oplus\Delta_m(W)=0,
\]
taking into account that $\Delta_{m-1}(W)=0$. Hence $S$ is an $m$-isometry in the case $m\geq3$.

For $m=2$, $S$ is defined in (\ref{eq:2.5}) with $W$ an isometry on $\mathcal{L}$, so $\Delta_1(S)=(T^*T-I+D)\oplus0$ on $\KK=\HH\oplus\mathcal{L}$, and further
\[
\Delta_2(S)=S^*\Delta_1(S)S-\Delta_1(S)
=(\Delta_2(T)+T^*DT-D)\oplus0=0,
\]
having in view the relation (\ref{eq:2.9}). So $S$ is a $2$-isometry.

It remains to prove that the lifting $S$ for $T$ is minimal, that meaning $\KK=\bigvee_{n\geq0}S^n\HH=:\KK_0$.

Indeed, let $k\in\KK$ be orthogonal to $\KK_0$. Since $k\perp\HH$ it follows that $k\in\mathcal{L}=\KK\ominus\HH$ and $k\perp S^n\HH$ for $n\geq1$. For $n=1$ we obtain from (\ref{eq:2.5}) that $k\perp E\HH$, hence $k\perp\overline{E\HH}=\mathcal{N}(W^*)$. Since $W$ is expansive, its range is closed, and therefore $\mathcal{N}(W^*)^\perp=\ran(W)$; consequently, $k\in\ran(W)$. Hence $k=W\ell_1$ for some element $\ell_1\in\mathcal{L}$. Next, since $k\perp S^2\HH$ we have for $h\in\HH$,
\[
0=(k,ET h\oplus WEh)=(W\ell_1,WEh)=(W^*W\ell_1,Eh).
\]
Therefore $W^*W\ell_1\perp\overline{E\HH}=\mathcal{N}(W^*)$, i.e. $W^*W\ell_1=W\ell'$ for some $\ell'\in\mathcal{L}$. Hence $\ell_1=(W^*W)^{-1}W\ell'\in\ran(W)$, because $\ran(W)=\{0\}\oplus\ls(\EE)$ (the weights $W_n$ are invertible) is a reducing subspace for the diagonal operator $W^*W$, and hence is invariant under $(W^*W)^{-1}$. It follows that there exists $\ell_2\in\mathcal{L}$ such that $\ell_1=W\ell_2$, so $k=W\ell_1=W^2\ell_2$.

By an induction argument one can show that $k\in\ran(W^n)$ for every integer $n\geq1$, hence $k=0$ ($W$ being a shift operator). We conclude $\KK=\KK_0$, which shows that $S$ is a minimal lifting for $T$. This completes the proof.
\end{proof}

Now we prove the identity (\ref{eq:2.15}) from the above proof.
\begin{lemma}\label{lem:identity}
Let $T\in\BH$ be as in Theorem \ref{thm:main}, and $D,D_n\in\BH$ be the operators from (\ref{eq:2.8}), respectively (\ref{eq:2.14}), for $n\geq0$. Then $D$ and $D_n$ satisfy the relation (\ref{eq:2.15}) for $m\geq3$ and every integer $n\geq1$.
\end{lemma}
\begin{proof}
For $n=1$ we have from the definition of $D_1$ in (\ref{eq:2.14}),
\[
\begin{aligned}
D_1&=T^*T+E^*E=\Delta_1(T)+I+D
=\binom10 I+\binom11\Delta_1(T)+D\\
&=\sum_{j=0}^{m-1}\binom1j\Delta_j(T)
+\left[1+\binom1{m-1}\right]D,
\end{aligned}
\]
taking into account that $E^*E=E_0^*E_0=D$ (by (\ref{eq:2.10})).

Now assume that the relation (\ref{eq:2.15}) holds for some $n\geq1$,
where $\binom{n}{p}=Q_p(n)$ has the meaning from the above proof (step 3). Then we first obtain
\begin{equation}\label{eq:2.16}
T^*D_nT=\sum_{j=0}^{m-1}\binom{n}{j}T^*\Delta_j(T)T
+\left[n+\binom{n}{m-1}\right]T^*DT.
\end{equation}
On the other hand, we have from (\ref{eq:2.14}) that $(D_nh,h)=\|S^n(h,0)\|^2$ for $h\in\HH$ and $n\geq1$. Since $E=JE_0$, where $J$ is the embedding into the zeroth coordinate of $\mathcal L=\ell_+^2(\EE)$, we have
\[
S(h,0)=(Th,(E_0h,0,0,\ldots)).
\]
Consequently,
\[
\begin{aligned}
(D_{n+1}h,h)&=\|S^{n+1}(h,0)\|^2
=\|S^n(Th,0)\|^2+\|W^n(E_0h,0,0,\ldots)\|^2\\
&=(D_nTh,Th)+(B_nE_0h,E_0h)
\end{aligned}
\]
for $h\in\HH$, where in the last term we used (\ref{eq:2.13}) for $j=0$. Since by (\ref{eq:2.10}) and (\ref{eq:2.11}), we also have
\[
E_0^*B_nE_0=(1+s_n)D+n\Delta_m(T)
=D+n\Delta_m(T)+\binom{n}{m-2}D,
\]
it follows from the above relations that
\[
D_{n+1}=T^*D_nT+\left[1+\binom{n}{m-2}\right]D+n\Delta_m(T).
\]
Next by using (\ref{eq:2.16}) as well as that $T^*\Delta_j(T)T=\Delta_j(T)+\Delta_{j+1}(T)$ and $T^*DT=D-\Delta_m(T)$ from (\ref{eq:2.9}), one obtains
\[
\begin{aligned}
D_{n+1}
&=\sum_{j=0}^{m-1}\binom{n}{j}\Delta_j(T)
+\sum_{j=0}^{m-1}\binom{n}{j}\Delta_{j+1}(T)\\
&\quad+\left[n+1+\binom{n}{m-1}+\binom{n}{m-2}\right]D
-\binom{n}{m-1}\Delta_m(T).
\end{aligned}
\]
But the second sum can be expressed as
\[
\sum_{j=0}^{m-1}\binom{n}{j}\Delta_{j+1}(T)
=\sum_{p=1}^{m}\binom{n}{p-1}\Delta_p(T)
=\sum_{p=1}^{m-1}\binom{n}{p-1}\Delta_p(T)
+\binom{n}{m-1}\Delta_m(T).
\]
Then we infer from the previous relation that
\begin{eqnarray*}
D_{n+1}
&=&\sum_{j=0}^{m-1}\binom{n}{j}\Delta_j(T)
+ \sum_{j=1}^{m-1}\binom{n}{j-1}\Delta_j(T)
+\left[n+1+\binom{n}{m-1}+\binom{n}{m-2}\right]D\\
&=&
I+\sum_{j=1}^{m-1}\binom{n+1}{j}\Delta_j(T)
+\left[n+1+\binom{n+1}{m-1}\right]D\\
&=&
\sum_{j=0}^{m-1}\binom{n+1}{j}\Delta_j(T)
+\left[n+1+\binom{n+1}{m-1}\right]D.
\end{eqnarray*}
So, by the induction argument we conclude that the relation (\ref{eq:2.15}) is true for every integer $n\geq1$, in the case $m\geq3$.
\end{proof}

\begin{remark}\label{rem:m2}
The case $m=2$ of Theorem \ref{thm:main} was already studied in \cite[Theorem~3.2]{BMS} where it was proved that every (2-)convex operator $T$ with the condition $\|T^n\|=O(\sqrt{n})$ as $n\to\infty$, has a $2$-isometric lifting $S$ with the ``invariant'' part $W$ an isometry. But by the above Theorem \ref{thm:main} one can choose $W$ to be a forward shift. Also, the proof in this case $m=2$ is different from that in \cite{BMS}, where the lifting $S$ was obtained by $n$-step successive convex liftings.
\end{remark}

Another interesting case in our context is $m=3$, and the corresponding result is mentioned in the following.

\begin{corollary}\label{cor:m3}
Let $T\in\BH$ be $3$-convex such that $\sup_{n\geq1}\|T^n\|/n<+\infty$. Then $T$ admits a minimal $3$-isometric lifting $S$ on a space $\KK=\HH\oplus\mathcal{L}$, such that $S|_{\mathcal{L}}$ is a $2$-isometric weighted shift on $\mathcal{L}$ and $S^{*n}S^n\HH\subset\HH$ for every integer $n\geq1$. Additionally, $S$ is expansive if $T$ is expansive.
\end{corollary}

This corollary shows that $W=S|_{\mathcal{L}}$ is just a $2$-isometry, but not an isometry in general. However, it is possible to get another lifting for $T$ with $W$ an isometry, for some $3$-convex operators $T$, as we see in the following.

\begin{example}\label{ex:shift}
Let $T$ on $\ls(\mathbb{C})$ be the weighted shift with the weights $\{w_n\}_{n\geq0}$, where $w_0=\sqrt6$ and
\[
w_n=\sqrt{\frac{p(n+1)}{p(n)}},\qquad
p(n)=n^2+n+4\qquad(n\geq1).
\]
It is easy to see that, if $\{e_n\}_{n\geq0}$ is the canonical orthonormal basis of $\ls(\mathbb{C})$, then for every integer $q\geq1$ we have
\[
\Delta_q(T)e_n=
\frac{\displaystyle\sum_{j=0}^{q}(-1)^{q-j}\binom{q}{j}p(n+j)}{p(n)}e_n
\qquad(n\geq0),
\]
taking into account that $p(0)=1$ because $w_0=\sqrt6$.

Clearly, $p(n)>0$ for $n\geq0$, $p(n+1)-p(n)=2(n+1)$ for $n\geq1$ and $p(1)-p(0)=5$. Thus, it follows that $(\Delta_1(T)e_0,e_0)=5$ and
\[
(\Delta_1(T)e_n,e_n)=\frac{2(n+1)}{p(n)}\qquad(n\geq1),
\]
therefore $\Delta_1(T)\geq0$, that is $T$ is expansive.

Concerning the operator $\Delta_2(T)$ we have $p(2)-2p(1)+p(0)=-1$ and
\[
p(n+2)-2p(n+1)+p(n)=2\qquad(n\geq1).
\]
So, for $\lambda=(\lambda_0,\lambda_1,\lambda_2,\ldots)$ with $\lambda_j\in\mathbb{C}$ one obtains
\[
\Delta_2(T)\lambda=\left(-\lambda_0,\frac{2}{p(1)}\lambda_1,
\frac{2}{p(2)}\lambda_2,\ldots\right).
\]
For $\Delta_3(T)$ we firstly have $p(3)-3p(2)+3p(1)-p(0)=3$ and
\[
p(n+3)-3p(n+2)+3p(n+1)-p(n)=0\qquad(n\geq1).
\]
This yields for $\lambda\in\ls(\mathbb{C})$ as above that
$\Delta_3(T)\lambda=(3\lambda_0,0,0,\ldots)$, which implies $\Delta_3(T)\geq0$. Hence $T$ is expansive and $3$-convex.

In addition, there exists a constant $\xi>0$ such that
\[
\|T^n\|^2=\sup_{q\geq0}\frac{p(n+q)}{p(q)}
\leq\xi(n+1)^2\qquad(n\geq0).
\]
So by Theorem \ref{thm:main}, $T$ admits an expansive $3$-isometric lifting with the ``orthogonal'' part of $T$ a $2$-isometry. But in this case we can get another lifting having this part even an isometry.

Indeed, for such purpose we firstly consider the operators $X,Y$ on $\ls(\mathbb{C})$ given by
\[
Xe_n=\frac{2n+5}{p(n)}e_n,\qquad
Ye_n=\frac{2}{p(n)}e_n\qquad(n\geq0).
\]
It is easy to see that
\begin{equation}\label{eq:2.17}
T^*XT-X=Y=T^*YT\geq0.
\end{equation}
Also, $X\geq\Delta_1(T)$ because $Xe_0=5e_0=\Delta_1e_0$, and
\[
(X-\Delta_1(T))e_n=\frac{3}{p(n)}e_n\qquad(n\geq1).
\]
Now consider $\HH_1=\ls(\mathbb{C})$ and $\KK_1=\HH_1 \oplus \ls(\HH_1)$. Let $V_1$ be the forward shift on $\ls(\HH_1)$, and $E_1:\HH_1\to \ls(\HH_1) $ be $E_1=J_0(X-\Delta_1(T))^{1/2}$, where $J_0:\HH_1\to \ls(\HH_1)$ is the embedding mapping. Then for the lifting $S_1$ of $T$ on $\KK_1$ with the representation
\[
S_1=\begin{pmatrix}T&0\\E_1&V_1\end{pmatrix},
\]
we have (by (\ref{eq:2.17})) the relations
\[
\Delta_1(S_1)=X\oplus0,\qquad
\Delta_2(S_1)=Y\oplus0,\qquad
\Delta_3(S_1)=0.
\]
Hence $S_1$ is a $3$-isometric lifting of $T$ with $S_1|_{\KK_1\ominus\HH_1}$ an isometry.
\end{example}

Notice that such liftings as $S_1$ were recently studied in \cite{CS,SuciuExpansive,SuciuTotoi}

\section{$m$-Isometric Liftings for $m$-Concave Operators}

In the sequel we show that a similar result to that of Theorem~\ref{thm:main} can be obtained for the $m$-concave operators, without any boundedness condition imposed. This extends the two-isometric lifting theory for concave operators in \cite{BS1,BS2} and is related to the expansive $m$-concave dilation result in \cite{Bu}. We begin with the following.

\begin{lemma}\label{lem:3.1}
If $T\in\BH$ is $m$-concave for an integer $m\geq2$, then $T$ is $(m-1)$-convex, that is $\Delta_m(T)\leq0$ implies $\Delta_{m-1}(T)\geq0$.
\end{lemma}
\begin{proof}
For $h\in\HH$ denote
\[
a_{n,p}:=(\Delta_p(T)T^nh,T^nh)
\]
for integers $n,p\geq0$, and briefly we write $a_n=a_{n,0}$.
Assume $\Delta_m(T)\leq0$, that is
\[
T^*\Delta_{m-1}(T)T-\Delta_{m-1}(T)=\Delta_m(T)\leq0.
\]
With the preceding notation this means
\[
a_{n,m}\leq0,\qquad a_{n+1,m-1}\leq a_{n,m-1}\qquad(n\geq0).
\]
In particular, $a_{n,m-1}\leq a_{0,m-1}=(\Delta_{m-1}(T)h,h)$, for $n\geq0$.

Assume (contrary to the conclusion) that $a_{0,m-1}<0$. Then using the relation $a_{n+1}-a_n=a_{n,1}$, as well as the identity
\[
\binom{n}{p}+\binom{n}{p-1}=\binom{n+1}{p},
\]
one obtains for $n\geq m-1$ that
\[
\begin{aligned}
a_n&=\sum_{j=0}^{m-2}\binom{n}{j}a_{0,j}
 +\sum_{p=0}^{n-m+1}\binom{n-p-1}{m-2}a_{p,m-1}\\
&\leq\sum_{j=0}^{m-2}\binom{n}{j}a_{0,j}
 +\sum_{p=0}^{n-m+1}\binom{n-p-1}{m-2}a_{0,m-1}\\
&=\sum_{j=0}^{m-2}\binom{n}{j}a_{0,j}
 +\binom{n}{m-1}a_{0,m-1}.
\end{aligned}
\]
Since $n\geq m-1$ we have $n^{-(m-1)}\sum_{j=0}^{m-2}\binom{n}{j}a_{0,j}\to0$ as $n\to\infty$. Therefore, for $n$ sufficiently large one obtains the contradiction
\[
0\leq \limsup_{n\to \infty} n^{-(m-1)}a_n\leq\frac{a_{0,m-1}}{(m-1)!}<0.
\]
Hence $(\Delta_{m-1}(T)h,h)=a_{0,m-1}\geq0$ for $h\in\HH$, that is $\Delta_{m-1}(T)\geq0$.
\end{proof}

Now we state and prove the main result of this section.

\begin{theorem}\label{thm:3.2}
Let $T\in\BH$ be $m$-concave for an integer $m\geq2$. Then $T$ admits a minimal $m$-isometric lifting $S_0$ on a Hilbert space $\KK_0=\HH\oplus\mathcal L_0$, of the form
\begin{equation}\label{eq:3.1}
S_0=\begin{pmatrix}T&0\\F&V\end{pmatrix},
\end{equation}
where $V$ is an expansive $m$-isometric weighted shift on $\mathcal L_0=\ell_+^2(\overline{\Delta_{m-1}(T)\HH})$, and $F:\HH\to\mathcal L_0$ is the operator
\begin{equation}\label{eq:3.2}
Fh=J_0[\Delta_{m-1}(T)]^{1/2}h\qquad(h\in\HH),
\end{equation}
while $J_0$ is the embedding mapping of $\overline{\Delta_{m-1}(T)\HH}$ into $\mathcal L_0$.
Furthermore, one has $S_0^{*n}S_0^n\HH\subset\HH$ for $n\geq1$, and $S_0$ is expansive if $T$ is expansive.
\end{theorem}
\begin{proof}
We use an argument similar to that in the proof of Theorem~\ref{thm:main}, in two stages. Assume $T\in\BH$ with $\Delta_m(T)\leq0$.

\emph{Step 1.} Define the operator $\Delta:\HH\to\mathcal D:=\overline{\Delta_{m-1}(T)\HH}$ by
\[
\Delta h=[\Delta_{m-1}(T)]^{1/2}h\qquad(h\in\HH),
\]
having in view that $\Delta_{m-1}(T)\geq0$ (by Lemma~\ref{lem:3.1}).

Since we have the relations
\[
0\leq-\Delta_m(T)=\Delta_{m-1}(T)-T^*\Delta_{m-1}(T)T
\leq\Delta_{m-1}(T)=\Delta^*\Delta,
\]
it follows that there exists a contraction $C\in\mathcal B(\mathcal D)$ such that
\[
C\Delta=(\Delta_{m-1}(T)-T^*\Delta_{m-1}(T)T)^{1/2}.
\]
Then the operator $A=C^*C$ satisfies the conditions
\begin{equation}\label{eq:3.3}
0\leq A\leq I_{\mathcal D},\qquad
\Delta^*A\Delta=\Delta_{m-1}(T)-T^*\Delta_{m-1}(T)T=-\Delta_m(T).
\end{equation}
Define the operator $F:\HH\to\mathcal L_0:=\ell_+^2(\mathcal D)$ in (\ref{eq:3.1}) by $F=J_0\Delta$, where $J_0d=(d,0,0,\ldots)\in\mathcal L_0$, for $d\in\mathcal D=\overline{\ran(\Delta)}$.

The operator $V$ in (\ref{eq:3.1}) will be a weighted shift with operator weights. To define these we consider the operator polynomials
\begin{equation}\label{eq:3.4}
A_n=I_{\mathcal D}+P(n)A,\qquad P(n)=n+\binom{n}{m-1}\qquad(n\geq0).
\end{equation}
Here $A$ is from (\ref{eq:3.3}), $A_0=I_{\mathcal D}\leq A_n$, so $A_n$ are expansive operators in $\mathcal B(\mathcal D)$ for $n\geq1$. Define the operator weights $V_n\in\mathcal B(\mathcal D)$ for $n\geq0$, and the unilateral shift $V\in\mathcal B(\mathcal L_0)$ by
\begin{equation}\label{eq:3.5}
V_n=A_{n+1}^{1/2}A_n^{-1/2},\qquad
V(d_0,d_1,\ldots)=(0,V_0d_0,V_1d_1,\ldots),\quad d_j\in\mathcal D.
\end{equation}
To show that the sequence $\{V_n\}_{n\geq0}$ is bounded, we first see that for $n\geq m-1$ the following relations hold:
\[
\begin{aligned}
P(n+1)&=n+1+\binom{n+1}{m-1}
=n+1+\frac{n+1}{n-m+2}\binom{n}{m-1}\\
&\leq n+1+m\binom{n}{m-1}
\leq m\left[n+\binom{n}{m-1}\right]=mP(n).
\end{aligned}
\]
For $1\leq n\leq m-2$ we have
\[
P(n+1)\leq n+2\leq m\leq mP(n).
\]
These together imply by (\ref{eq:3.4}) that
\[
A_n\leq A_{n+1}\leq(m+1)A_n\qquad(n\geq1).
\]
Thus, for $V_n$ from (\ref{eq:3.5}) we get the relations
\[
I\leq V_n^*V_n=A_n^{-1/2}A_{n+1}A_n^{-1/2}\leq(m+1)I\qquad(n\geq0),
\]
having also in view that $A_0=I$ and $V_0^*V_0=A_1\leq I+2A\leq3I$. Hence the sequence $\{V_n\}_{n\geq0}$ is bounded, and $\|V\|\leq\sqrt{m+1}$.

In addition, because $V_n^*V_n\geq I$ for $n\geq0$, it follows that $V^*V\geq I$, that is $V$ is an expansive weighted shift. It remains to show that $V$ is an $m$-isometry. For this purpose denote by $J_q:\mathcal D\to\mathcal L_0$ the natural inclusion of $\mathcal D$ into $q$-coordinate of $\mathcal L_0$, for $q\geq0$. Since $\ran(J_q)\perp\ran(J_{q'})$ for $q\neq q'$, the operator $V^{*n}V^n$ has a diagonal block matrix on $\mathcal L_0=\bigoplus_{q\geq0}\mathcal D$ for $n\geq0$, and
\begin{equation}\label{eq:3.6}
J_q^*V^{*n}V^nJ_q=A_q^{-1/2}A_{n+q}A_q^{-1/2}\qquad(n,q\geq0).
\end{equation}
Now for each fixed $q$, $A_{n+q}$ is an operator polynomial of degree at most $m-1$ (from (\ref{eq:3.4})), therefore
\begin{equation}\label{eq:3.7}
\sum_{n=0}^{m}(-1)^{m-n}\binom{m}{n}A_{n+q}=0\qquad(q\geq0).
\end{equation}
But this ensures by (\ref{eq:3.6}) and the structure of $V$ that $\Delta_m(V)=0$, that is $V$ is an $m$-isometry.

Also notice that $J_0^*V^{*n}V^nJ_0=A_n$ for $n\geq0$, while the diagonal structure of $V^{*n}V^n$ gives $J_q^*V^{*n}V^nJ_0=0$ for $q\geq1$ and $n\geq0$. Because $J_0\mathcal D=\mathcal{N}(V^*)$ we obtain the conditions
\begin{equation}\label{eq:3.8}
V^{*n}V^n\mathcal{N}(V^*)\subset\mathcal{N}(V^*)\qquad(n\geq1).
\end{equation}
Finally, having in view that $V$ is a forward shift on $\mathcal L_0$, it follows that $\bigcap_{n\geq0}\ran(V^n)=\{0\}$, meaning that $V$ is analytic.

\emph{Step 2.} Define the operator $S_0$ on the space $\KK_0=\HH\oplus\mathcal L_0$ as in (\ref{eq:3.1}) with entries $F$ and $V$ above. Clearly, $S_0$ is a lifting of $T$, and using the representation (\ref{eq:3.1}) we obtain for $n\geq1$,
\begin{equation}\label{eq:3.9}
S_0^n=\begin{pmatrix}T^n&0\\F_n&V^n\end{pmatrix},\qquad
S_0^{*n}S_0^n=\begin{pmatrix}G_n&0\\0&V^{*n}V^n\end{pmatrix},
\end{equation}
where
\begin{equation}\label{eq:3.10}
G_n=T^{*n}T^n+F_n^*F_n,\qquad
F_n=\sum_{j=0}^{n-1}V^{n-j-1}FT^j.
\end{equation}
Here we take into account that $V^{*n}F_n=0$, by using that $\overline{\ran(F)}=J_0\mathcal D=\mathcal{N}(V^*)$ and the conditions (\ref{eq:3.8}). So we have $S_0^{*n}S_0^n\HH\subset\HH$ for $n\geq1$.

To get an expression for $G_n$ similar to that of $D_n$ in (\ref{eq:2.15}), we compare the first equality in (\ref{eq:3.3}) satisfied for $\Delta_{m-1}(T)\geq0$ and $\Delta_m(T)\leq0$, with the relation (\ref{eq:2.9}) for $D\geq\Delta_m(T)\geq0$. Thus, by using the relations 
$$
F^*V^{*n}V^nF=\Delta_{m-1}(T)- P(n) \Delta_m (T)
$$
and 
$$
G_{n+1}=T^*G_nT +\Delta_{m-1}(T)-P(n)\Delta_m (T),
$$
one can obtain (as in Lemma~\ref{lem:identity}) for $G_n$ the expression 
\[
G_n=\sum_{j=0}^{m-2}\binom{n}{j}\Delta_j(T)
+\left[n+\binom{n}{m-1}\right]\Delta_{m-1}(T)\qquad(n\geq1).
\]

This shows that $G_n$ is a polynomial in $n$ of degree at most $m-1$, hence
\[
\sum_{n=0}^{m}(-1)^{m-n}\binom{m}{n}G_n=0.
\]
Since $\Delta_m(V)=0$, it follows (by (\ref{eq:3.9})) that $\Delta_m(S_0)=0$, that is $S_0$ is an $m$-isometry on $\KK_0$. Also, as in the proof of Theorem~\ref{thm:main} one can see that $S_0$ is a minimal lifting for $T$, taking into account that $V$ is analytic. In addition, if $T^*T\geq I$ then $G_1\geq I$ by (\ref{eq:3.10}) and so $S_0^*S_0\geq I$ by (\ref{eq:3.9}), i.e. $S_0$ is expansive. This completes the proof.
\end{proof}

\begin{remark}\label{rem:3.3}
The case $m=2$ meaning that $T$ is (2-) concave was studied in \cite{BS1,BS2}, where one proves that $T$ has a $2$-isometric lifting $S_0$ with $V=S|_{\KK_0\ominus\HH}$ an isometry.

In the case $m\geq 3$, a result corresponding to Theorem~\ref{thm:3.2} was recently obtained in \cite[Theorem~3.1]{Bu}, under the additional condition that $T$ is expansive, which is omitted in the case $m=3$. This appears to be a special case of Theorem~\ref{thm:3.2}, but the above proof of this theorem is different from that in \cite{Bu}.

In addition, we remark that if $\Delta_m(T)\neq0$ in Theorem~\ref{thm:3.2}, then $A\neq0$ in (\ref{eq:3.3}). Since $J_0^*V^{*n}V^nJ_0=A_n$ for $n\geq0$, we obtain (by (\ref{eq:3.7})), for $m\ge 3$ that
\[
\begin{aligned}
J_0^*\Delta_{m-1}(V)J_0
&=\sum_{n=0}^{m-1}(-1)^{m-n-1}\binom{m-1}{n}J_0^*V^{*n}V^nJ_0\\
&=\sum_{n=0}^{m-1}(-1)^{m-n-1}\binom{m-1}{n}A_n=A\neq0.
\end{aligned}
\]
Hence $\Delta_{m-1}(V)\neq0$, meaning that $V$ (and so $S_0$) is a strict $m$-isometry.
\end{remark}

\begin{remark}\label{rem:3.4}
From the relations (\ref{eq:3.10}) we have $\|T^n\|\leq\|S_0^n\|$ for $n\geq0$. As $S_0$ is an $m$-isometry it follows that $T$ satisfies the boundedness condition (\ref{eq:2.4}) inferred from $S_0$. Hence this condition is necessary for an $m$-concave operator $T$, and this does not need to be imposed in the hypothesis of Theorem~\ref{thm:3.2}.

We know (by Lemma~\ref{lem:3.1}) that if $\Delta_m(T)\leq0$ then $\Delta_{m-1}(T)\geq0$. But the $m$-convexity does not imply the $(m-1)$-concavity. For example, for $T=2I$ we have $\Delta_m(T)=3^mI$ for every $m\geq1$. However, under a boundedness condition, this converse implication can occur, and the above theorems assure the corresponding liftings. 
\end{remark}
\bigskip

\begin{proposition}\label{prop:3.5}
Let $T\in\BH$ be $m$-convex for an integer $m\geq3$. The following statements hold:
\begin{itemize}
 \item[(i)] If $n^{-(m-1)}\|T^nh\|^2\xrightarrow[n\to\infty]{}0$ for every $h\in\HH$, then $T$ is $(m-1)$-concave and it admits an $(m-1)$-isometric lifting.

 \item[(ii)] If $\sup_{n\geq1}n^{-(m-2)}\|T^n\|^2<+\infty$, then $T$ is $(m-2)$-convex and it admits an $(m-1)$-isometric lifting.
\end{itemize}
\end{proposition}
\begin{proof}
Fix $h\in\HH$ and, as in the proof of Lemma~\ref{lem:3.1}, put
\[
a_{n,p}:=(\Delta_p(T)T^nh,T^nh),\qquad n,p\geq0.
\]
Since $T$ is $m$-convex, we have
\[
a_{n+1,m-1}-a_{n,m-1}
=(\Delta_m(T)T^nh,T^nh)\geq0,
\]
so the sequence $\{a_{n,m-1}\}_{n\geq0}$ is nondecreasing. The finite-difference identity used in the proof of Lemma~\ref{lem:3.1} gives, for $n\geq m-1$,
\[
a_{n,0}=\sum_{j=0}^{m-2}\binom{n}{j}a_{0,j}
+\sum_{p=0}^{n-m+1}\binom{n-p-1}{m-2}a_{p,m-1}.
\]
Next, since $a_{p,m-1}\geq a_{0,m-1}$ and by using the identity for binomial coefficients from the proof of Lemma \ref{lem:3.1}, we get
\[
a_{n,0}\geq\sum_{j=0}^{m-2}\binom{n}{j}a_{0,j}
+\binom{n}{m-1}a_{0,m-1}.
\]

Assume that $a_{0,m-1}>0$. Dividing the preceding inequality by $n^{m-1}$ and letting $n\to\infty$, we obtain
\[
0=\lim_{n\to\infty}\frac{a_{n,0}}{n^{m-1}}
\geq\lim_{n\to\infty}\left[
\sum_{j=0}^{m-2}\frac{\binom{n}{j}}{n^{m-1}}a_{0,j}
+\frac{\binom{n}{m-1}}{n^{m-1}}a_{0,m-1}
\right]
=\frac{a_{0,m-1}}{(m-1)!}>0,
\]
which is a contradiction. Here we used the hypothesis in (i) and the equality $a_{n,0}=\|T^nh\|^2$. Thus $a_{0,m-1}\leq0$ for every $h\in\HH$, and hence $\Delta_{m-1}(T)\leq0$. Theorem~\ref{thm:3.2}, applied with $m-1$ in place of $m$, now yields an $(m-1)$-isometric lifting of $T$. This proves (i).

For (ii), let
\[
M:=\sup_{n\geq1}n^{-(m-2)}\|T^n\|^2<+\infty.
\]
Then, for every $h\in\HH$ and $n\geq1$,
\[
0\leq n^{-(m-1)}\|T^nh\|^2
\leq\frac{M\|h\|^2}{n}\xrightarrow[n\to\infty]{}0.
\]
Therefore the hypothesis of (i) is satisfied. It follows that $T$ is $(m-1)$-concave and admits an $(m-1)$-isometric lifting. Finally, applying Lemma~\ref{lem:3.1} with $m-1$ in place of $m$, we obtain that $T$ is $(m-2)$-convex. This proves (ii) and completes the proof.
\end{proof}

This proposition suggests an alternation between convexity and concavity, in some conditions. More precisely, we infer the following.

\begin{corollary}\label{cor:3.6}
Let $T\in\BH$ be an operator satisfying for the integers $p,m$ with $2\leq p < m$,
\begin{equation}\label{eq:3.11}
\sup_{n\geq1}\frac{1}{n^{p-1}}\|T^n\|^2<+\infty.
\end{equation}
 The following statements hold:

(i) If $T$ is $m$-convex, then $(-1)^{m-q}\Delta_q(T)\geq0$ for $p\leq q < m$. 

(ii) If $T$ is $m$-concave, then $(-1)^{m-q}\Delta_q(T)\leq0$ for $p\leq q < m$. 

Furthermore, if $T$ satisfies the hypothesis of (i) or (ii), then $T$ admits a $q$-isometric lifting, for every integer $q$ such that $p\le q < m$, and this lifting is expansive if $T$ is such.
\end{corollary}

\begin{proof}
The assertions (i) and (ii) result by applying successively Proposition~\ref{prop:3.5} and Lemma~\ref{lem:3.1}, taking into account from (\ref{eq:3.11}) that $n^{-q}\|T^n\|^2\to0$ as $n\to\infty$ for $p\leq q < m$. 

The assertion of (i) shows that if $\Delta_m(T)\geq0$ then $\Delta_q(T)\geq0$, when $m-q$ is even, and $\Delta_q (T) \le 0$ when $m-q$ is odd. In the first case the condition (\ref{eq:3.11}) ensures, by Theorem~\ref{thm:main}, that $T$ admits a $q$-isometric lifting, which is even expansive if $T$ is such.

In the case that $\Delta_q(T)\leq0$ we have by Theorem~\ref{thm:3.2} that $T$ admits a $q$-isometric lifting which is expansive when $T$ is such.

A similar argument can be used for the last assertion, when $\Delta_m(T)\leq0$ (as in (ii)). Here $\Delta_q(T)\leq0$ when $m-q$ is even, and $\Delta_q(T)\geq0$ if $m-q$ is odd. We omit other details.
\end{proof}

\begin{remark}\label{rem:3.7}
The condition (\ref{eq:3.11}) is essential for the above assertions (i) and (ii). In this sense we have seen that the operator $T$ from Example~\ref{ex:shift} is expansive (i.e. $1$-convex) and $3$-convex, but it is neither $2$-convex nor $2$-concave. In this case Corollary~\ref{cor:3.6} (i) cannot be applied, because the condition (\ref{eq:3.11}) is not satisfied for $m=3$ in the case $p=q=2$, taking in view that $\sup_{n\geq1}\frac{1}{\sqrt{n}}\|T^n\|=+\infty$, for $T$ from Example~\ref{ex:shift}.

Now we give two examples of operators which show that the convexity and the concavity cannot be transferred from an order to the other. 
\end{remark}

\begin{example}\label{ex:3.8}
Let $T$ on $\ell_+^2(\mathbb C)$ be the weighted shift with the weights $\{w_n\}_{n\geq0}$, where
\[
w_n=\sqrt{\frac{h(n+1)}{h(n)}},\qquad
h(n)=(n+1)^2+\frac{(-1)^n}{4}\qquad(n\geq0).
\]
We have $h(n)>0$ for $n\geq0$, and for the orthonormal basis $\{e_n\}_{n\geq0}$ in $\ell_+^2(\mathbb C)$,
\[
\Delta_q(T)e_n=
\frac{\displaystyle\sum_{j=0}^{q}(-1)^{q-j}\binom{q}{j}h(n+j)}{h(n)}e_n
\qquad(q\geq1,\ n\geq0).
\]
It is easy to see that $h(n+1)-h(n)=2n+3-\frac{(-1)^n}{2}>0$ and
\[
h(n+2)-2h(n+1)+h(n)=2+(-1)^n>0,
\]
\[
h(n+3)-3h(n+2)+3h(n+1)-h(n)=-2(-1)^n,
\]
for $n\geq0$. These show us, by the above relations of $\Delta_q(T)$, that $\Delta_1(T)\geq0$, $\Delta_2(T)\geq0$, so $T$ is expansive and $2$-convex. But it is neither $3$-convex nor $3$-concave. Hence the expansivity and 2-convexity imply neither 3-convexity nor 3-concavity. 
\end{example}

\begin{example}\label{ex:3.9}
Let $T$ on $\ell_+^2(\mathbb C)$ be the weighted shift having the weights
\[
t_n=\sqrt{\frac{q(n+1)}{q(n)}},\qquad
q(n)=n+2-3^{-n}+\frac18\left(-\frac13\right)^n\qquad(n\geq0).
\]
Then one can show that for $n\geq0$ and $p\geq2$ we have $q(n)>0$ and
\[
b_{p,n}:=\sum_{j=0}^{p}(-1)^{p-j}\binom{p}{j}q(n+j)
=(-2)^p3^{-(n+p)}\bigl(-1+2^{p-3}(-1)^n\bigr).
\]
It follows that $b_{2,n}<0$ for $n\geq0$, which yields $\Delta_2(T)\leq0$, so $T$ is $2$-concave. But $b_{3,n}=0$ for $n$ even and $b_{3,n}>0$ if $n$ is odd, which implies that $\Delta_3(T)\geq0$, i.e. $T$ is $3$-convex.

In the case $p=4$ we see that $b_{4,n}>0$ if $n$ is even and $b_{4,n}<0$ when $n$ is odd, so $T$ is neither $4$-convex, nor $4$-concave. 

Since $T$ is $3$-convex and $2$-concave, the first assertion of (i) in Corollary~\ref{cor:3.6} and the condition (\ref{eq:3.11}) are verified in the case $m=3$ and $p=q=2$. Thus in this case one can apply Corollary~\ref{cor:3.6} (i).
In addition we see that neither 4-convexity nor 4-concavity is assured by 2-concavity and 3-convexity, even if a certain boundedness condition of type \eqref{eq:3.11} of order $p=2$ holds.
\end{example}

\bigskip

{\small
	{\em Author's address}: 
	
	{\em Laurian Suciu}, Department of Mathematics and Informatics,
	
	 ``Lucian Blaga'' University of Sibiu, Sibiu, Romania.
	
	e-mail: \texttt{laurians2002@yahoo.com}.
}


\begin{thebibliography}{99}

\bibitem{AbdullahLe}
B. Abdullah and T. Le,
\emph{The structure of $m$-isometric weighted shift operators},
Oper. Matrices \textbf{10} (2016), no.~2, 319--334.
DOI: 10.7153/oam-10-17.

\bibitem{Agler}
J. Agler,
\emph{A disconjugacy theorem for Toeplitz operators},
Amer. J. Math. \textbf{112} (1990), no.~1, 1--14.
DOI: 10.2307/2374849.

\bibitem{AglerStankusI}
J. Agler and M. Stankus,
\emph{$m$-isometric transformations of Hilbert space, I},
Integral Equations Operator Theory \textbf{21} (1995), no.~4, 383--429.
DOI: 10.1007/BF01222016.

\bibitem{AglerStankusII}
J. Agler and M. Stankus,
\emph{$m$-isometric transformations of Hilbert space, II},
Integral Equations Operator Theory \textbf{23} (1995), no.~1, 1--48.
DOI: 10.1007/BF01261201.

\bibitem{AglerStankusIII}
J. Agler and M. Stankus,
\emph{$m$-isometric transformations of Hilbert space, III},
Integral Equations Operator Theory \textbf{24} (1996), no.~4, 379--421.
DOI: 10.1007/BF01191619.

\bibitem{BMS}
C. Badea, V. M\"uller and L. Suciu,
\emph{High order isometric liftings and dilations},
Studia Math. \textbf{258} (2021), no.~1, 87--101.
DOI: 10.4064/sm200330-25-8.

\bibitem{BS1}
C. Badea and L. Suciu,
\emph{The Cauchy dual and $2$-isometric liftings of concave operators},
J. Math. Anal. Appl. \textbf{472} (2019), no.~2, 1458--1474.
DOI: 10.1016/j.jmaa.2018.12.002.

\bibitem{BS2}
C. Badea and L. Suciu,
\emph{Hilbert space operators with two-isometric dilations},
J. Operator Theory \textbf{86} (2021), no.~1, 93--123.
DOI: 10.7900/jot.2020feb05.2298.

\bibitem{BermudezMartinonNegrin}
T. Berm\'udez, A. Martin\'on and E. Negr\'in,
\emph{Weighted shift operators which are $m$-isometries},
Integral Equations Operator Theory \textbf{68} (2010), no.~3, 301--312.
DOI: 10.1007/s00020-010-1801-z.

\bibitem{Bu}
M. Bucha\l{}a,
\emph{Every expansive $m$-concave operator has $m$-isometric dilation},
Linear Algebra Appl. \textbf{732} (2026), 93--107.
DOI: 10.1016/j.laa.2025.12.001.

\bibitem{BuchalaShift}
M. Bucha\l{}a,
\emph{$m$-isometric weighted shifts with operator weights},
Bull. Malays. Math. Sci. Soc. \textbf{48} (2025), Article 71.
DOI: 10.1007/s40840-025-01853-1.

\bibitem{CS} A. Crăciunescu and L. Suciu, \emph{Brownian extensions in the context of three-isometries}, J. Math. Anal. Appl. 529 (2024), no. 1, Article 127591, 19 pp.
DOI: https://doi.org/10.1016/j.jmaa.2023.127591

\bibitem{MajdakBrownianConcave}
W. Majdak and L. Suciu,
\emph{Brownian isometric parts of concave operators},
New York J. Math. \textbf{25} (2019), 1067--1090.

\bibitem{MajdakBrownianParts}
W. Majdak and L. Suciu,
\emph{Brownian type parts of operators in Hilbert spaces},
Results Math. \textbf{75} (2020), Article 5, 1--23.
DOI: 10.1007/s00025-019-1130-8.

\bibitem{MajdakSuciuConvex}
W. Majdak and L. Suciu,
\emph{Convex and expansive liftings close to two-isometries and power bounded operators},
Linear Algebra Appl. \textbf{617} (2021), 1--26.
DOI: 10.1016/j.laa.2021.01.009.

\bibitem{MajdakSuciuTriangulations}
W. Majdak and L. Suciu,
\emph{Triangulations of operators with two-isometric liftings},
Integral Equations Operator Theory \textbf{93} (2021), Article 10, 1--24.
DOI: 10.1007/s00020-021-02625-9.

\bibitem{SuciuTwoLiftings}
L. Suciu,
\emph{On operators with two-isometric liftings},
Complex Anal. Oper. Theory \textbf{14} (2020), Article 5, 1--16.
DOI: 10.1007/s11785-019-00960-9.

\bibitem{SuciuExpansive}
L. Suciu,
\emph{Operators with expansive $m$-isometric liftings},
Monatsh. Math. \textbf{198} (2022), 165--187.
DOI: 10.1007/s00605-021-01648-z.

\bibitem{SuciuTotoi}
L. Suciu and E. A. Totoi,
\emph{Three-isometric liftings with invariant isometric part},
Banach J. Math. Anal. \textbf{15} (2021), Article 66.
DOI: 10.1007/s43037-021-00150-w.

\bibitem{SFb} B. Sz.-Nagy, C. Foias, H. Bercovici, L. K\'erchy, \emph{Harmonic Analysis of Operators on Hilbert Space},
2nd ed., Revised and enlarged edition. Universitext, Springer, New York, 2010. Zbl 1234.47001, 
\url{https://doi.org/10.1007/978-1-4419-6094-8}

\end{thebibliography}
\end{document}